\documentclass[12pt,emtex]{book}
\usepackage{times}

\usepackage{graphicx}
\usepackage{amsfonts}       %
\usepackage{amsmath} 
\usepackage{amssymb}        %
\usepackage{amsthm}         

\usepackage{epsfig}
\psfigdriver{dvips}

\usepackage{latexsym}
\usepackage[matrix,curve]{xypic}

\usepackage[german]{babel}  
\newcommand{\A}{\mathbb{A}}
\newcommand{\N}{\mathbb{N}}
\newcommand{\Z}{\mathbb{Z}}
\newcommand{\Q}{\mathbb{Q}}

\newcommand{\C}{\mathbb{C}}
\newcommand{\F}{\mathbb{F}}

\begin{document}
\centerline{\large \bf Siegel modular forms mod $p$}

\bigskip\noindent
\centerline{R.Weissauer}

\thispagestyle{empty}
\bigskip\noindent

\bigskip\noindent

\bigskip\noindent
Let $R$ be a commutative ring with unit. Assume $1/2\in R$. For a
prime $p$ let $\F_p$ denote the prime field of characteristic $p$
and $\overline\F_p$ an algebraic closure. We always assume $p\neq
2$.

\bigskip\noindent
{\it The moduli stack $A_{g,R}$}. Let $g\geq 2$ be an integer.
Consider the moduli stack $A_{g,R}$  of isomorphism classes of the
following data $(S,A,\lambda)$
\begin{itemize} \item[]  a ring homomorphism $R\to S$ of commutative rings, \item[]  an abelian scheme $A\to
Spec(S)$ of relative dimension $g$ together with its zero section
$e:Spec(S)\to A$, \item[]  a principal polarization $\lambda: A\to
A^\vee$ over $Spec(S)$.
\end{itemize}
To $(S,A,\lambda)$ we associate the sheaf of invariant relative
differential forms $\Omega_A=e^*(\Omega^1_{A/S})$, which is a
projective rank 1 module $\omega = \Lambda^g(\Omega_A)$ on $S$
functorial with respect to base change.

\bigskip\noindent
{\it Modular forms}. By definition a modular form of integral
weight $k$ and genus $g$ over $R$ is  an assignment
$$ (S,A,\lambda) \mapsto f(S,A,\lambda) \in \Gamma(Spec(S),\omega^{\otimes k}) $$
with the following functorial properties \begin{itemize}
\item[(i)] Basechange compatibility: For ring homomorphisms of commutative rings $S\to S'$
one has
$$ f(S',A\times_{Spec(S)} Spec(S'),\lambda \times_{Spec(S)} id) =  f(S,A,\lambda) \otimes_{{\cal O}_S} 1 \ .$$
\item[(ii)] For $S$-isomorphic data $(A,\lambda)$ and
$(A',\lambda')$ one has
$$ f(S,A,\lambda) = f(S,A',\lambda') \ .$$
\end{itemize}
Let  $M_{g,R}^k$ be the space of modular forms of weight $k$. The
direct sum over all integers $k$ defines the graded $R$-algebra
$M_{g,R}^\bullet = \bigoplus_k M_{g,R}^k$ of modular forms.

\bigskip\noindent
{\it Variants}. One has variants for the definition of modular
forms: a) additional level-N-structure over extension rings $R$ of
$\Z[\frac{1}{N},\zeta_N]$. This defines the stack $A_{g,R,N}$. For
this we assume $2p\nmid N$ later in case $R$ is the spectrum of a
field $\kappa$ of characteristic $p$. Or b) tensor valued modular
forms, i.e. sections $f(S,A,\lambda) \in
\Gamma(Spec(S),\Omega_A^{\otimes k})$ with corresponding
functorial properties.

\bigskip\noindent
{\it The classical case}. $ M_{g,\C}^k$ is isomorphic to the space
of all holomorphic functions $f:H_g \to \C$ on the Siegel upper
half space $H_g$ of genus $g$ with the property $$ f(\gamma
\langle\Omega\rangle) = det(C\Omega +D)^k f(\Omega) $$ for all
$\gamma \in (\begin{smallmatrix} A & B\cr C &
D\cr\end{smallmatrix}) \in \Gamma_g= Sp_{2g}(\Z)
 $. The isomorphism is given by
$$ f \mapsto f(\Omega) = f(\C,\C/\Gamma_\Omega,\lambda_{can}) $$
for periods $\Omega$ in the Siegel upper half space $H_g$ and the
standard polarization $\lambda_{can}$ defined for the standard
period lattice $\Gamma_\Omega=\Z^g+\Omega\Z^g$.

\bigskip\noindent
{\it The Hasse invariant for $R=\F_p$}. Let $Lie(A/S)\cong
\Omega_A^*$ denote the Lie algebra of translation invariant
relative vectorfields $X$, and let denote $D=D_X$ the derivation
associated to $X$. Consider the Hasse-Witt map
$$ \Phi: Lie(A/S) \to Lie(A/S) $$
$$ \Phi(D) = D^p  \ . $$ Notice $\Phi$ is $p$-linear, i.e.
$\Phi(\lambda v)=\lambda^p \Phi(v)$ for $\lambda \in S$. For an
arbitrary $p$-linear endomorphism $\Phi:V\to V$ of a locally free
$S$-module $V$ with $pS=0$ put $V^{(p)} = V\otimes_S S$ such that
$v\lambda\otimes_S \lambda' = v \otimes_S \lambda^p\lambda'$. Then
$\Phi$ extends to a $S$-linear morphism $\Phi: V^{(p)} \to V$. Let
$ST^p(V)$ denote the $p$-symmetric elements in the tensor product.
It contains the symmetrization $Sym^p(\otimes^p V)$ as a
$S$-submodule. The $S$-linear morphism $v\otimes\lambda \mapsto
\lambda \cdot v\otimes v \otimes \cdots \otimes v$ ($p$ copies)
induces an $S$-isomorphism $V^{(p)}
\overset{\sim}{\longrightarrow} ST^p(V)/Sym^p(\otimes^p V)$ by
[D], I \S 10. If we compose the projection $$ST^p(V)\to
ST^p(V)/Sym^p(\otimes^p V)$$ with the inverse of this isomorphism
and then with $\Phi$, we obtain a $S$-linear morphism $ST^p(V)\to
V$. It can be considered as an element in $$V \otimes_S S^p(V^*)\
,$$ where $S^p(V^*)$ is the symmetric tensor quotient of the
$p$-th tensor product of $V^*$ defined as the quotient by the
equivalence relation induced from the symmetry relations
$...\otimes x\otimes ... \otimes y \otimes .. \sim ...\otimes
y\otimes ... \otimes x \otimes .. $.

\bigskip\noindent
Applied for the Hasse-Witt map $\Phi$ this construction yields a
sheaf homomorphism $ ST^p(Lie(A/S)) \to Lie(A/S) $ or equivalently
an element
$$ \A \in \Gamma(Spec(S),S^p(\Omega_A) \otimes \Omega_A^*) \ .$$
Exterior powers of $p$-linear maps are $p$-linear. Hence the
exterior powers $\Lambda^i(\A)$ of $\A$ define vector valued
modular forms in characteristic $p$. In particular for $i=g$ this
defines the Hasse invariant $A$ as a modular form of weight
$k=p-1$ in the space $M_{g,\F_p}^{p-1}$ defined above
$$ A \in \Gamma(Spec(S),\omega^{p-1}) \ .$$

\bigskip\noindent
{\it The ungraded ring of modular forms}. Assume $\kappa$ to be
$\F_p$ or its algebraic closure $\overline\F_p$. Let denote
$M_{g,\kappa}$ the associated ungraded ring of modular forms over
$\kappa$ obtained from the graded ring $M_{g,\kappa}^\bullet$.
This ring is a normal domain.  Then

\bigskip\noindent
{\bf Unique Faktorization Theorem}. {\it There exists a prime
$p_0$ such that for $p>p_0$ the normal domain $M_{g,\kappa}$ of
modular forms over $\kappa$ is a unique factorization ring for all
$g\geq 2$}.

\bigskip\noindent
The lemma is related to the following graded version for modular
forms over $\kappa$ of some level $N\geq 1$, where we tacitly
assume a root of unity $\zeta_N$  to be chosen in $\kappa$.

\bigskip\noindent
{\bf Theorem}. {\it For $p>p_0$ and $g\geq 3$ and $N\geq 1$ and
$2p\nmid N$ the Picard group of the moduli stack $A_{g,\kappa,N}$
defined over $\kappa$ is
$$ Pic(A_{g,\kappa,N}) \cong \Z \cdot c_1(\omega)\ \oplus\
sp_{2g}(N\Z/N^2\Z) \ .$$}

\medskip\noindent {\it Remark}. In particular this group is of rank
one, and generated by the Chern class $c_1(\omega)$ modulo a
finite torsion subgroup isomorphic to the \lq Lie algebra\rq
$sp_{2g}(\Z/N\Z)$ of the symplectic group over the ring $\Z/N\Z$.
The corresponding assertion ist false for $g=2$. It does not hold
for $g \geq 3$ and levels of type $N=2^i$ either. However there
are similar statements in this latter case.

\bigskip\noindent
{\it Remark}. To obtain the results above on Picard groups for
$\F_p$ it is enough to prove this over $\overline \F_p$ using
descend via Hilbert theorem 90 and the vanishing of the Brauer
groups of these fields.

\bigskip\noindent
{\it Fourier expansion}. The $q$-expansion of a modular form is
obtained by evaluating forms at the Mumford families of
principally polarized abelian varieties defined by torus
quotients. The $q$-expansion principle [CF] gives injective maps
$$ M_{g,R}^k \hookrightarrow R[[q]]:=R[[q_{ij}]] \quad , \quad q_{ij}=q_{ji},
1\leq i,j\leq g $$ defined by
$$ f \mapsto \sum_T a_f(T) q^T $$
with Fourier coefficients $a_f(T)\in R$ attached to half-integral
symmetric matrices $T$ (written as symbolic exponents). Notice,
the relevant Mumford family defines an abelian scheme over
$R[[q]][\frac{1}{q_{ij}}]$. Nevertheless by Koecher's effect $f\in
R[[q]]$.

\bigskip\noindent
The theory of Chai-Faltings [CF] implies that $M_{g,\Z}^k = \{
f\in M_{g,\C}^k \ \vert \ a_f(T)\in \Z\}$. This relates modular
form over $\Z$ to classical Siegel modular forms.

\bigskip\noindent
{\it Example}. For $R=\F_p$ the Hasse invariant $A\in
M_{g,\F_p}^{p-1}$ has Fourier expansion $A=1$. In other words
$a_f(T)=0$ for all $T\neq 0$.

\bigskip\noindent
We will see below that the last theorem implies

\bigskip\noindent
{\bf Proposition}. {\it Assume $p> p_0$ for $p_0$ as in the lemma.
Then the kernel of the ring homomorphism
$$ \varphi: M_{g,\F_p} \to \F_p[[q]] $$
of ungraded rings induced from the $q$-expansion is the principal
ideal generated by $(A-1)$ for the Hasse invariant $A$.}

\bigskip\noindent
Let $S_{g,\F_p} \subseteq \F_p[[q]] $ denote the image of the
homomorphism $\varphi$. It is not hard to see for $p\geq 2g+2$
that $S_{g,\F_p}$ coincides with the ring obtained by reduction
mod $p$ from the ring of the integral $q$-expansion of all modular
forms in $M_{g,\Z}$. Furthermore one has the following rather
obvious properties: The subring $S_{g,\F_p}$ of the ring of power
series $ \F_p[[q]] $ is stable under the following operators
\goodbreak
\begin{enumerate}
\item[(1)] Hecke operators $T(l)$ for primes $l$ different from $p$
\item [(2)] The operators
$$ U: \sum_T a(T)q^T\ \mapsto\ \sum_T a(pT)q^T $$
$$  V: \sum_T a(T)q^T\ \mapsto\ \sum_T a(T)q^{pT} \ .$$
\item [(3)]
$$\Phi: \F_p[q_{11},..,q_{gg}]]\ \mapsto\
\F_p[[q_{11},..,q_{g-1g-1}]] $$
$$ q_{*g}=q_{g*}\mapsto 0 $$
and $q_{ij}\mapsto q_{ij}$ for $i\neq g,j\neq g$. Notice
$$\Phi(S_{g,\F_p})=S_{g-1,\F_p}$$ and $$\Phi(M_{g,\F_p})\ \subseteq\
M_{g-1,\F_p}\ .$$ \end{enumerate} For the operator $V$ use that
reduction modulo $p$ of the normalized Hecke operator $T(p)$
coincides with $V$ for weights $k\geq g+1$. Of course one can
assume $k\geq g+1$ without restriction of generality, since
multiplication with the Hasse invariant $A$ commutes with the
Hecke operators $T(l)$ for all $(l,p)=1$. For forms $f$ with
coefficients in $\F_p$ we have $f\vert V = f^p$.

\bigskip\noindent
{\it Proof of the proposition}. For the elliptic case $g=1$ see
[Sw]. The general case is reduced to the elliptic case as follows.
First notice
$$ dim(S_{g,\F_p})\ \geq\ dim(M_{g,\F_p}) - 1 \ .$$
In fact, if $\kappa$ is a field and $B=\bigoplus_{i=0}^\infty B_i$
is a finitely generated integral (graded) $\kappa$-algebra of
Krull dimension $dim(B)$, and if $\varphi: B\to S$ is a ring
homomorphism onto an integral domain $S$ such that $\varphi\vert
B_i$ is injective for every $i$, then $dim(S)\geq dim(B) -1$.
Since the stack $A_{g,\F_p}$ is irreducible by a result of Chai,
the graded ring of modular forms satisfies the assumptions made on
$B$. Hence this general fact from commutative algebra can be
applied for the graded ring of modular forms. Obviously
$Ker(\varphi)$ is a prime ideal. Therefore it is a prime ideal of
height one because the kernel is nontrivial, since $0\neq A-1\in
Ker(\varphi)$. Hence $Ker(\varphi)$ is a principal ideal by the
unique factorization theorem stated above.

\bigskip\noindent
Since $A-1$ is contained in the kernel of $\varphi$, it is enough to show that $A-1$ is
irreducible. If not $A-1 =FG$ for $F=F_0+\cdots +F_s$ and $G=G_0+\cdots +G_t$ and $p-1=r+s,
r>0, s>0$ and forms $F_i,G_i$ of weight $i$. But then $\Phi^{g-1}(A-1)=A-1$ (in genus $g=1$)
implies $\Phi^{g-1}(F_s)\neq 0$ and $\Phi^{g-1}(G_t)\neq 0$ by degree reasons (notice the Hasse
invariant $A$ is nonvanishing of weight $p-1$ for genus $g=1$). But this factorization of $A-1$
for genus $g=1$ would contradict the result of Swinnerton-Dyer [Sw]. This proves the
proposition. \hfill $\square$

\bigskip\noindent
By a similar argument, using the $\Phi$-operator to reduce to the
case $g=2$, we have the following

\bigskip\noindent
{\bf Corollary}. {\it For $g\geq 2$ the Hasse invariant is an
irreducible modular form. In particular the locus of singular
abelian varieties (the zero locus of the Hasse invariant $A$) is
irreducible.}

\bigskip\noindent
For $g=2$ this follows from results of Koblitz.

\goodbreak
\bigskip\noindent
{\bf Corollary}. {\it Suppose $g\geq 2$, $p>p_0$ and $f\in
M_{g,\F_p}^k$. Then the following statements are equivalent
\begin{enumerate}
\item[1)] $f$ is totally $p$-singular, i.e. $f=\sum_T a(T) q^T\ $ and $a(T)=0\ $ unless $p\vert T$.
\item[2)] $f=A^r h^p$ for some $r\in \N_0$ and some $h\in
M_{g,\F_p}^{k'}$. \end{enumerate}}

\bigskip\noindent
{\it In particular $f=0$ for $k<p-1$, if one of these two
equivalent conditions holds.}

\bigskip\noindent
{\it Proof}. To show that 1) implies 2) we may assume, that $A$ does not divide $f$. If $f$ has
weight $k$, we choose $l$ such that $K=k+l(p-1)\geq g+1$. Then $A^lf$ has weight $K\geq g+1$,
hence also $(A^lf)\vert U$. Then $((A^lf)\vert U)\vert V= ((A^lf)\vert U)^p$ has weight $pK$,
and it has the same Fourier expansion as $f$. Therefore $((A^lf)\vert U)^p = A^m f$ for
$m=k+lp$ by the $q$-expansion principle, since both sides have the same weight and the same
$q$-expansion. So at least the power $A^{lp}$ divides $((A^lf)\vert U)^p$. Hence $A^l$ divides
$(A^lf)\vert U$ and $h_1=A^{-l}(A^lf)\vert U)$ is a modular form of weight $k$ such that $h_1^p
= A^kf$. Therefore the irreducible form $A$ divides $h_1$. In other words $h_1=Ah_2$, hence
$A^ph_2^p = A^k f$. Since $A$ does not divide $f$ by assumption, we conclude $h_2^p = A^{k-p}
f$ and either $k=p$, or again $h_2=Ah_3$ holds for some modular form $h_2$. If we continue,
this must terminate after finitely many steps. Hence $k$ is divisible by $p$ and there exists a
modular form $h=h_\nu$ such that $h^p=f$. \hfill $\square$

\bigskip\noindent
{\it Question}. Is $S_{g,\F_p}\subseteq \F_q[[q]]$ preserved by
the differential operator $\sum_T a(T)\cdot q^T \mapsto \sum_T
det(T)a(T)\cdot q^T$ ? Similarly is $\sum_T T\cdot a(T) q^T$ a
vector valued modular form?

\bigskip\noindent
{\it Proof of the theorem}. For the proof we can extend the base
field to become $\kappa=\overline \F_p$. Use Hilbert theorem 90
and vanishing of the Brauer group for finite fields, and later for
fields of fractions of Henselian  discrete valuation rings with
algebraically closed residue field. We can also replace $N$ by a
suitably large integer not divisible by $2p$, so that $X$ is a
smooth variety without restriction of generality.

\bigskip\noindent
1) {\it Change of level}. Assume we are over a base scheme $S$
over $\Z[\frac{1}{N},\zeta_N]$ for $N\geq 3$. Consider the natural
morphism $A_{g,S,N} \to A_{g,S}$, which is equivariant with
respect to the group $G_N=Sp(2g,\Z/N\Z)$. Pullback defines a
natural homomorphism $Pic(A_{g,S}) \to Pic(A_{g,S,N})$. A line
bundle on $A_{g,S}$ is given by its pullback ${\cal L}$ on
$A_{g,S,N}$ together with descent data $\varphi_\sigma:
\sigma^*({\cal L}) \cong {\cal L}$ for $\sigma\in G_N$ satisfying
the obvious cocycle conditions. Since $H^0(A_{g,S,N},{\cal O}^*)=
{\cal O}^*_S(S)$ and the group $G_N$ is perfect, the cocycle datum
is uniquely determined, once it exists. Hence there exists an
exact sequence
$$ 0 \to Pic(A_{g,S}) \to Pic(A_{g,S,N})^{G_N} \to H^2(G_N,{\cal
O}^*_S(S)) $$ For odd $N$ and $g\geq 3$ we claim $H^2(G_N,C)=0$
for $C={\cal O}^*_S(S)$. We may replace $C$ by an arbitrary finite
trivial $G_N$-module. For perfect groups $G$ the K\"{u}nneth theorem
and the universal coefficient theorem imply $H^2(G,C) =
Hom(H_2(G,\Z),C)$. So it is enough to show vanishing of the Schur
multiplier $M(G)=H^2(G,\Q/\Z)\cong H_2(G,\Z)^D$. Since the Schur
multiplier is additive for products of perfect groups, one can
reduce to the case where $N$ is a prime power. According to the
group Atlas, see also[FP], this Schur multiplier vanishes in the
case where $N$ is a prime different from two. Using results of
[FP] one reduces this to show $H^2(G_N,\Z/l\Z)$ vanishes. For
$(l,N)=1$ see [FP]. So we can assume $C=\F_p$. For $N=p^r$ for
$r\geq 2$ and $p$ odd now we proceed by induction on $r$ using
$G_{p^r}/Lie =: G = G_{p^{r-1}}$ with respect to the abelian group
$Lie=sp_{2g}(\Z/p\Z)$. The Hochschild-Serre spectral sequence and
the induction assumption together with $H^0(G,H^2(Lie))=0$ and
$H^1(G,Lie^D)=0$ then prove the induction step. For $p\neq 2$ one
has an exact sequence $0\to Lie^D \to H^2(Lie,\F_p)\to
\Lambda^2(Lie^D)\to 0$, hence the first statement
$H^0(G,H^2(Lie))=0$ follows from the obvious facts $(Lie^D)^G=0$
and $\Lambda^2(Lie^D)^G=0$. The second assertion $H^1(G,Lie^D)=0$
is less obvious. It follows by induction on the genus $g$ using
restriction to parabolic subgroups and standard techniques of
group cohomology. Hence under our assumptions on $N$ the pullback
defines an isomorphism
$$ \fbox{$ Pic(A_{g,S}) \cong Pic(A_{g,S,N})^{G_N} $} \ .$$

\bigskip\noindent
By the methods above one also obtains the vanishing of cohomology
groups of the finite group $G_N=Sp_{2g}(\Z/N\Z)$ for odd $N$ in
the following special cases
$$ H^i\bigl(G_N,\C^*\bigr) =0 \quad , \quad g\geq 3, \ i=1,2 $$
$$ H^i\bigl(G_N,(sp_{2g}(\Z/N\Z),Ad)\bigr) =0 \quad , \quad g\geq 3, \ i=0,1 \mbox{ and } \ 2\nmid N $$
$$ H^3\bigl(G_N,\Z)_{l-torsion} =0 \quad , \quad \mbox{for } g\geq 3 \mbox{ and } \ 2\nmid l \ .$$
In the proof these vanishing theorems can be used to allow the
passage from level $N=1$ to level $N>1$ and vice versa.

\bigskip\noindent
2) {\it  Analytic results over $\C$}.  For $g\geq 3$ we have a)
For $\Gamma_g=Sp_{2g}(\Z)$
$$ Pic(\Gamma_{g}) :=H^1(H_g/\Gamma_g,{\cal O}^*(H_g)) \cong
H^2(\Gamma_g,\Z) = \Z\cdot c_1(\omega)\ ,$$  and b) For all
principal congruence subgroups $\Gamma=\Gamma_g(N)$ for $N\geq 3$
$$ H^2(\Gamma,\Z) \cong Pic_{an}(H_g/\Gamma) \cong
Pic(A_{g,\C,N})\ ,$$ and c) $$Pic(\Gamma)_{tors} \cong
(\Gamma^{ab})^D\ .$$ See [F], p. 257 and Hilfssatz 4. Finally d),
that $Pic(\Gamma)/\Z c_1(\omega)$ is a finite group for $g\geq 3$
and arbitrary subgroups $\Gamma$ of finite index in $\Gamma_g$.
This follows from theorems of A.Borel and vanishing  theorems
(resp. $L^2$-vanishing  theorems for $g=3$) for Lie algebra
cohomology.

\bigskip\noindent
Recall $G_N=\Gamma_g/\Gamma_g(N)$. By step 1 the $E_2$-term of the
Hochschild-Serre spectral sequence $ H^i(G_N, H^j(\Gamma_g(N),\Z))
\Longrightarrow H^{i+j}(\Gamma_g)$ for $i+j=2$ has only $H^0(G_N,
H^2(\Gamma_g(N),\Z))$ as nontrivial summand, and degenerates in
this degree  at the $E_3$-term, i.e. the differential $d_3$
induces a canonical exact sequence $$ 0\to H^2(\Gamma_g,\Z) \to
H^2(\Gamma_g(N),\Z)^{\Gamma_g} \to H^3(G_N,\Z)\ .$$ Since $
H^3(G_N,\Z) \cong H^2(G_N,\Q/\Z)$ the vanishing of Schur
multiplier explained in step 1 implies $ H^3(G_N,\Z) =0$ for odd
$N$ and $g\geq 3$. Hence the restriction map induces a canonical
isomorphism $$ \fbox{$ res: \ H^2(\Gamma_g,\Z) \ \cong\
H^2(\Gamma_g(N),\Z)^{\Gamma_g}$}
$$ for $g\geq 3$ and odd integers $N$. Hence both for the stack setting
as well as for  the analytic setting the descend from level $N$ to
level one only amounts to consider $G_N$-invariants. This is used
in step 4) below.

\bigskip\noindent
A first application of this descend principle: For $N\geq 3$ and
$g\geq 3$ the long exact sequence of $G_N$-invariants for the
exact sequence $0 \to ({\Gamma^{ab}})^D \to Pic(A_{g,\C,N})
\overset{p}{\to} \Z \to 0 $ from assertion 2c) and 2d) and the
vanishing of the cohomology groups $H^i(G_N,({\Gamma^{ab}})^D)=0$
for $i=0,1$ implies that $p$ induces an isomorphism $ p:
Pic(A_{g,\C,N}^{G_N}) \cong \Z $. Here we used ${\Gamma_g(N)^{ab}}
= \Gamma_g(N)/\Gamma_g(2N) \cong sp_{2g}(\Z/N\Z)$ and the
vanishing theorems of step 1, in particular that the perfect group
$G_N$ acts trivially on the cyclic quotient group $\Z$. The
descend principle above thus gives the canonical isomorphism
$$ p: Pic(A_{g,\C}) \cong \Z \ .$$ Since $c_1(\omega)$ generates
$Pic(A_{g,\C}) \cong H^2(\Gamma_g,\Z)$ by 2b) we conclude that the
projection map $p$ induces an isomorphism
$$ Pic(A_{g,\C,N})\ \cong\ \Z \cdot c_1(\omega) \oplus Pic(A_{g,\C,N})_{tors}
 $$ for odd $N$ and $g\geq 3$. We would like to carry this over
 to the Picard groups of $A_{g,\F_p,N}$ via
smooth and proper base change for etale $l$-adic cohomology for
all $l\neq p$. For this consider $N\geq 3$ and compactifications.

\bigskip\noindent
3) {\it Passage to compactification}. There exists an exact
$G_N$-equivariant sequence
$$ 0\to \Z^{r} \to Pic(\overline{A_{g,\kappa,N}}) \to
Pic(A_{g,\kappa,N}) \to 0 $$ for suitable toroidal
compactification $\overline{A_{g,\kappa,N}}$ of $A_{g,\kappa,N}$
in the sense of Chai-Faltings. Indeed the subgroup generated by
the boundary divisors is torsion-free by Koecher's principle,
hence isomorphic to $\Z^r$ for a suitable integer $r$ depending on
$N,g$ and the compactification. The same sequences exist over
$\Z[\frac{1}{N},\zeta_N]$ and $\C$.

\bigskip\noindent
4) {\it Base change}. The strategy of proof is to reduce the
assertion of the theorem to the case of level $N=1$ using step 1).
For $N=1$ we have to show $Pic(A_{g,\kappa})=\Z\cdot c_1(\omega)$.
For this choose some auxiliary prime $l$ different from $p$, and
first show $Pic(A_{g,\kappa})\otimes \Z_l=\Z_l\cdot c_1(\omega)$.
If the stack $A_g$ were a smooth and proper scheme, this could be
easily deduced from the analytic facts stated in 2) by proper and
smooth base change. Since it is not one has to pass to some larger
level $N$ and consider a toroidal compactification for which it is
possible to apply the smooth and proper base change theorem. Using
step 1) and 3) one can pass back to derive
$Pic(A_{g,\kappa})\otimes \Z_l=\Z_l\cdot c_1(\omega)$. For varying
$l$ this implies that $Pic(A_{g,\kappa})/(\Z\cdot c_1(\omega))$ is
a finite $p$-primary abelian group.

\bigskip\noindent
To make this outline more precise choose an integer $l$ prime to
$p$. One can assume $N\geq 3$. First one replaces $\F_p$ by its
algebraic closure $\overline\F_p$. Apply the snake lemma to the
exact sequences of step 3 with respect to multiplication by $l$.
This gives a diagram
$$ \xymatrix{0\ar[r] & Pic(\overline X/\kappa)[l] \ar[d] \ar[r] &
 Pic(A_{g,\kappa,N})[l]\ar[d] \ar[r] & (\Z/l\Z)^r
 \ar@{=}[d] \ar[r] & Pic(\overline  X/\kappa)/l \ar[d] \cr
0\ar[r] & Pic(\overline X/\C)[l] \ar[r] &
 Pic(A_{g,\C,N})[l] \ar[r] & (\Z/l\Z)^r \ar[r] & Pic(\overline X/\C)/l   } $$
 for toroidal compactifications $\overline  X$ as in step 3.
Notice $Pic(\overline  X/\kappa)\cong H^1(\overline
 X/\kappa,\mu_l)$ for $\kappa=\overline\F_p$ and $\kappa=\C$. Hence
the left vertical arrow is an isomorphism by the proper smooth
etale base change theorem. This implies that the second vertical
specialization map is injective. Next notice that the abelian
variety $Pic^0(\overline  X/\kappa)$ vanishes. Look at $m$-torsion
points for large auxiliary integers $m$ prime to $p$. Then
$Pic^0(\overline  X/\kappa)$ is a group of type $(\Z/m\Z)^{2d}$.
As a subgroup of $Pic(\overline  X/\kappa)$ this forces $d=0$. To
show this use base change for the first etale cohomology to reduce
this to the case $\kappa =\C$, where this immediately follows from
step 2). The vanishing of $Pic^0(\overline  X/\kappa)$ implies
$Pic(\overline  X/\kappa)=NS(\overline  X/\kappa)$. By the Kummer
sequence one has inclusions $Pic(\overline  X/\kappa)/l
\hookrightarrow H^2_{et}(\overline  X/\kappa,\mu_l)$ respectively
canonical inclusions $Pic(\overline  X/\overline\F_p)\otimes\Z_l
\hookrightarrow H^2_{et}(\overline  X/\overline\C,\Z_l)$ for prime
$l$. By base change $H^2_{et}(\overline
 X/\overline\F_p,\mu_l)\cong H^2_{et}(\overline  X/\C,\mu_l)$ respectively
$H^2_{et}(\overline  X/\overline\F_p,\Z_l)\cong H^2_{et}(\overline
 X/\C,\Z_l)$. Furthermore $H^2_{et}(\overline  X/\C,\Z_l) \cong
H^2_{et}(\overline  X/\C,\Z)\otimes\Z_l$. Hence the smooth and
proper base change theorem implies that the fourth vertical map of
the diagram above is injective. This implies by the 5-lemma, that
the second vertical map also is surjective. Hence the second
vertical map is an isomorphism
$$ Pic(A_{g,\overline\F_p,N})[l] \ \cong\ Pic(A_{g,\C,N})[l] \ .$$
The argument also gives an injective map
$$ Pic(A_{g,\overline\F_p,N})/l \ \hookrightarrow\ Pic(A_{g,\C,N})/l \ .$$
Together this implies the analog of assertion 2d) over
$\overline\F_p$. Hence
$$Pic(A_{g,\overline\F_p,N})/Pic(A_{g,\overline\F_p,N})_{tors}
\cong \Z\ .$$ Since $Pic(A_{g,\overline\F_p,N})_{tors}$ is
$(\Gamma_g(N)^{ab})^D$ up to $p$-power torsion as shown above, by
taking $G_N$-invariants the argument at the end of step 2 give the
following result
$$ Pic(A_{g,\C,N})\otimes\Z_l\ \cong\ (\Z_l \cdot c_1(\omega))\ \oplus\ (\Z_l\otimes
Pic(A_{g,\C,N})_{tors}) \ . $$ for all primes $l$ different from
$p$. Warning: Usually taking invariants does not commutes with the
tensor product with $\Z_l$. Here it does by a simple inspection.

\bigskip\noindent
5) {\it $p$-index}. To control the cotorsion of the subgroup
generated by $c_1(\omega)$ in the group $Pic(A_{g,\F_p,N})$ for
$N=1$ one suitably embeds $A_{2,\F_p} \hookrightarrow A_{g,\F_p}$
such that $\omega$ pulls back to $\omega$. Suppose ${\cal L}^e
=\omega$ for a line bundle ${\cal L}$ on $A_{g,\F_p}$. Then of
course the same holds for the pullback of ${\cal L}$ to
$A_{2,\F_p}$. So this reduces us to consider the case $g=2$ for
this question. Now Igusa [I] has given an explicit description of
the coarse moduli space $A_2^{coarse}$ of the stack $A_2$. In fact
for $p\neq 2,3$ the space $A_{2,\F_p}^{coarse}$ is a quotient
$U/G$, where $U$ is the open complement of some closed subset of
codimension 2 in the projective space ${\mathbb P}^3$, and where
$G$ is the finite group $G=\Z/2\Z\times \Z/3\Z \times \Z/5\Z
\times \Z/6\Z$ acting on $U$. The same holds over
$\kappa=\overline\F_p$. Consider
$$ \xymatrix{ & A_{2,\overline\F_p,3}\ar[dl]\ar[dr] & \cr A_{2,\overline\F_p} & & A_{2,\overline\F_p}^{coarse} } \ .$$
A suitable power of $\omega^r$ descends to
$A_{2,\overline\F_p}^{coarse}$. Since the Picard group of $U$ is
cyclic, this provides an a priori bound $e_0$ for the integer $e$
independent from $g$ and $p$. Since on the other hand $e$ must be
a power of $p$ by step 4) we find that for large enough $p>p_0$
necessarily $e=1$ holds. This proves $Pic(A_{g,\overline\F_p}) =
\Z c_1(\omega) \oplus T$ for a finite $p$-primary torsion group
$T$.

\bigskip\noindent
6) {\it $p$-torsion}. To show that the torsion $T$ vanishes it is
enough to control the $p$-torsion $Pic(A_{g,\kappa})[p]$. For this
we may pass to a suitably large level $N$ not divisible by $2p$.
Then the covering space $X=A_{g,\kappa,N}$ is smooth. It is enough
to show $Pic(X)[p]=0$. In general, for smooth $X$ over $\kappa$
one can control $Pic(X)[p]=0$ by the closed 1-forms on $X$. In
fact, if $C$ denotes the Cartier operator, then one has
 a surjection
$$ \bigg\{ \eta\in H^0(X,\Omega_X^1)\ \vert\ d\eta=0,C\eta=\eta\bigg\}
\cong H^1\bigl(X_{flat},\mu_p\bigr)\ \ \twoheadrightarrow \ \
Pic(X)[p] \ .$$ Hence to show that $Pic(X)[p]$ vanishes under the
assumptions made, it is enough to proof the following

\bigskip\noindent
{\bf Vanishing Theorem}. {\it Suppose $g\geq 2$ and $p\neq 2$ and
let $\kappa$ be $\F_p[\zeta_N]$ or $\overline\F_p$. Then for
$X=A_{g,\kappa,N}$ and arbitrary level $N$ not divisible by $p$
the following holds
$$ H^0(X,\Omega_X^1) =0 \ .$$}

\bigskip\noindent
{\it Remark}. The proof of this vanishing theorem can be extended
to the case of regular alternating differential forms $\Omega_X^i$
on $X$ of degree $i$ for $i< \frac{g}{2}$ and $i<\frac{p-1}{2}$.

\bigskip\noindent
{\it Remark}. We could compute $p_0$ from step 5 of the above
proof of the theorem and the assumption $p\neq 2$.

\bigskip\noindent
{\it Proof of the vanishing theorem}. The proof proceeds by
induction on $g$. In some sense the case $g=2$ is the most
complicated case. For the proof we can choose some embedding
$\Z[\zeta_N]$ to $\overline\F_p$ and extend the base field
$\kappa$ to become $\overline \F_p$. We can also replace $N$ by a
suitably large integer not divisible by $p$, so that $X$ is a
smooth variety without restriction of generality.

\bigskip\noindent
a) {\it $q$-expansion}. As a first step in the proof use that any
regular 1-form $\eta$ on $A_{g,\F_p,N}$ automatically extends to a
1-form $\eta\in \Omega_{\overline X}^1(log)$ with log-poles on the
toroidal compactification $\overline X =
\overline{A_{g,\kappa,N}}$ (Koecher's principle). Furthermore
notice the Kodaira-Spencer isomorphism in the sense of [CF], p.107
$$ \Omega_{\overline X}^1(log) \cong  S^2(\Omega_A) \ .$$
In other words, one can expand $\eta$ in the following form using
$q$-expansion
$$ \eta\ =\ \sum_{T\geq 0} Trace\bigl(a(T)\cdot dlog(q)\bigr)\cdot q^T \ ,$$
where $a(T)$ is a symmetric g by g matrix with coefficients in
$\kappa$ and where $T$ runs over the halfintegral semidefinite g
by g matrices (notation is chosen to hide $N$). E.g. for $g=2$
$$ dlog(q) = \begin{pmatrix}
\frac{dq_{11}}{q_{11}} & \frac{dq_{12}}{q_{12}} \cr
\frac{dq_{21}}{q_{21}} & \frac{dq_{22}}{q_{22}} \cr \end{pmatrix}
\ .$$

\bigskip\noindent
b) {\it Equivariance}. For unimodular $U\in Gl(g,\Z)$ congruent to
$id$ modulo $N$ with image $\overline U$ in $Gl(g,\Z/p\Z)$ the
Fourier coefficients $a(T)$ of the 1-form $\eta$ satisfy
$$ (*) \quad \quad \fbox{$ a({}^t U T U) = {}^t \overline U a(T) \overline
U $} \ .$$ In particular $a(0)=0$.

\bigskip\noindent
c) {\it Cartier operator}. Suppose there exists a 1-form $\eta\neq
0$ on $X$. Without restriction of generality one can choose $\eta$
such that there exists $T$ with $a(T)\neq 0$ and $p\nmid T$.
Notice otherwise the form is closed, so one can apply the Cartier
operator $C: \Omega_{X,closed}^1 \to \Omega_{X^{(p)}}^1$. Since
$C(\sum_T a(T)q^T dlog(q)) = \sum_T \tilde a(T)q^{T/p} dlog(q)$,
this process terminates after finitely many steps and gives the
desired form $\eta$ on $X^{(p^s)}$. We may replace $X^{(p^s)}$ by
$X$, since this amounts just to another choice of
$\zeta_N\in\kappa$. By b) we can in addition assume $T =
(\begin{smallmatrix}
* & * \cr * & \nu \end{smallmatrix})$ with last entry $T_{gg}=\nu \not\equiv
0$ mod $p$.

\bigskip\noindent
d) {\it A boundary chart}. For the sake of simplicity temporarily
assume $N=1$. Let $$\pi: {\cal A}_{g-1}\to A_{g-1}$$ be the
universal abelian scheme with polarization $\lambda: {\cal
A}_{g-1} \to {\cal A}_{g-1}^\vee$. On the pullback
$\Xi=\Xi_\xi=j^*({\cal P}^{-1})$ of the Poincare ${\mathbb
G}_m$-torsor ${\cal P}$ over ${\cal A}_{g-1}\times_{A_{g-1}} {\cal
A}_{g-1}^\vee$ under $j=id\times \lambda$ the group ${\mathbb
G}_m$ acts such that the canonical map $u:\Xi\to {\cal A}_{g-1}$
induces an isomorphism $\Xi/{\mathbb G}_m\cong {\cal A}_{g-1}$. A
formal boundary chart $S=S_\xi$ of the toroidal compactification
is now obtained by the formal completion of $\overline \Xi =
(\Xi\times \A^1)/{\mathbb G}_m$ along the zero section. Indeed
$$u: \overline \Xi\to {\cal A}_{g-1}$$ is a line bundle over
${\cal A}_{g-1}$. The complement of the zero section is the
${\mathbb G}_m$-bundle $\Xi=(\Xi\times {\mathbb G}_m)/{\mathbb
G}_m$ we started from. See [CH], p.104 ff for $X_\xi=\Z$, $E_\xi=\
{\mathbb G}_m$ and $E_\xi=\A^1$.

\bigskip\noindent
e) {\it Fourier-Jacobi expansion}. Then $\eta$, considered as a
differential form on $S$ with log-poles, can be expanded $
u_*(\eta\vert_S) \ =\ \sum_{\nu} u_*(\eta)^{\chi_\nu} q^\nu $ with
respect to the characters $\chi_\nu(t)=t^\nu $ of ${\mathbb G}_m$.
Here $q\in \Gamma(\overline \Xi,{\cal L}_0^{-1})$ is the
tautological section for ${\cal L}_0=(id\times\lambda)^*({\cal
P})$. Notice $q=\tau(\chi_1,\chi_1)$ in the notation of [CF],
p.106. To relate this to the Fourier-Jacobi expansion one needs to
consider the semiabelian principally polarized scheme
${}^\heartsuit G$ defined over $S$ by Mumford's construction. It
contains a split torus $T$ over $S$ of rank one with the abelian
scheme quotient ${}^\heartsuit G/T\cong {\cal A}_{g-1}$. Notice
$\Omega_T\cong {\cal O}_S$ via $\frac{dq}{q} \mapsto 1 $. The
exact sequence $0\to \Omega_{{\cal A}_{g-1}} \to
\Omega_{{}^\heartsuit G} \to \Omega_T \to 0$ of relative
differential forms over $S$ and the Kodaira-Spencer isomorphism
$Symm^2(\Omega_{{}^\heartsuit G}) \cong \Omega^1(log)$ give rise
to an isomorphism $\Omega^1_S = u^*({\cal E})$ for a vector bundle
${\cal E}$ on ${\cal A}_{g-1}$, which is pinched into an exact
sequence
$$ 0 \to \Omega^1_{{\cal A}_{g-1}} \to {\cal E} \overset{r}{\to} \ {\cal O}_{{\cal
A}_{g-1}} \to 0 \ .$$ The coefficients $u_*(\eta)^\chi$ are global
sections of $u_*(u^*({\cal E}))^\chi = {\cal E}\otimes_{{\cal
O}_{{\cal A}_{g-1}}} u_*({\cal O}_S)^\chi$.

\bigskip\noindent
{\it Separability}. Notice $u_*({\cal O}_S)^\chi = {\cal
L}_0^{\nu}$, where ${\cal L}_0=(id\times\lambda)^*({\cal P})$ is a
totally symmetric line bundle on ${\cal A}_{g-1}$ normalized along
the zero section. Recall $p\neq 2$. Hence for $$\nu\not\equiv 0
\mod p
$$ the line bundle ${\cal L}={\cal L}_0^\nu$ is
separable over the base $A_{g-1,\kappa}$ of characteristic $p$.
This allows to relate the sections of $u_*({\cal O})^\chi$ to
theta functions. This will be explained in the next section.
Indeed on some etale extension $S$ of the base (e.g.
$S=A_{g-1,\kappa,2\nu}$) the line bundle ${\cal L}$ is of type
$\delta=(2\nu,..,2\nu)$ in the sense of Mumford (see also the next
section). By further increasing the level we will see that then
$({\cal A}_{g-1},\pi,{\cal L})$ admits a $\delta$-marking in the
following sense.

\bigskip\noindent
f) {\it Theta functions}. For  an abelian scheme ${\cal A}\to S$ of relative dimension $g$ and
a  totally symmetric relative ample line bundle ${\cal L}$ on ${\cal A}$ normalized along the
zero section Mumford  defined the theta group scheme ${\cal G}({\cal L})$ over $S$ ([M2], prop
1). For a type $\delta=(d_1,..,d_g)$ of even integers with $d_i\vert d_{i+1}$ for all $i$ one
also has the Heisenberg group scheme ${\cal G}(\delta)$ over $S$ as defined in [M2], p.77.
Here, as in [M2], we will suppose a separability condition, namely that $d=\prod_{i=1}^g d_i$
is invertible on the base scheme $S$. Then by definition a symmetric $\vartheta$-structure, or
$\delta$-marking, for $({\cal A},\pi,{\cal L})$ is an isomorphism ${\cal G}({\cal L}) \cong
{\cal G}(\delta)$ with some additional compatibility properties as defined in [M2], p.79. Let
$V(\delta)$ be the free $\Z[d^{-1}]$-module generated of functions $K(\delta)=\bigoplus_i
(\Z/d_i\Z)\to \Z[d^{-1}]$. Then [M2] prop. 2 and the discussion in loc. cit p.81 implies for a
given $\delta$-marking that there exists an equivariant isomorphism of ${\cal O}_S$-sheaves
$$ \alpha: V_\delta \otimes_\Z K \ \cong \
 \pi_*({\cal L}) $$
 for some invertible sheaf $K$ on $S$, unique up to multiplication
 with some unit in ${\cal O}^*(S)$. Furthermore [CF], theorem 5.1
 gives $K^{8d^4}\otimes \omega^{4d^4}\cong {\cal O}_S$, since $det(
\pi_*({\cal L}))=K^d$ for $d= (2\nu)^g$. Hence $K^2\otimes\omega$
is a torsion line bundle on $S$ killed by $4d^4$. We notice, that
Mumford's theory of theta functions also gives an injections
$$ K \hookrightarrow {\cal O}_S \ .$$

\bigskip\noindent
{\it Additional remark}. Suppose $g\geq 3$. By suitably extending the level $N$ to large $N'$,
but still with $(N',p)=1$, all torsion line bundles in $Pic(S)[4d^4]$ become trivial over
$A_{g-1,\kappa,N'}$. For $g\geq 3$ this is first shown over $\C$, but then follows by the
methods used for the proof of the theorem. Alternative argument (for all $g\geq 1$). Use that
$K$ is defined over $\Z[\frac{1}{2\nu},\zeta_{N'}]$. Since the characteristic $p$ fiber of
$A_{g,\kappa,N'}$ is irreducible for $(N',p)=1$, hence a principal divisor, one can reduce the
statement to characteristic zero, where this can be explicitly computed by cocycles.

\bigskip\noindent
g) {\it Existence of $\delta$-markings}. In the case of base
fields existence of $\delta$-markings for a given symplectic
isomorphism
 $g_2:H({\cal L}^2) \cong H(2\delta)$ is shown
in [M] \S 2 (remark 3, p. 319). Such a symplectic isomorphism
$g_2$ exists over $S=A_{g-1,\kappa,N}$ by extending the level $N$.
It is enough that $4\nu$ divides the level $N$. For a general base
$S$, different to the case of a base field considered in loc.
cit., there exists an obstruction to lifting $g_2$ to a symmetric
isomorphism $f_2:{\cal G}({\cal L}^2)\to {\cal G}(2\delta)$. The
proof of loc. cit. carries over, if this obstruction vanishes. The
obstruction is defined by the map
$$\delta: {\cal H}({\cal L}^2) \to  Pic(S)[8\nu] $$
given by $\delta(s)= t_{s^*}({\cal L}^2)\otimes {\cal L}^{-2}$ for
sections $s:S\to {\cal H}({\cal L}^2)$. This map is quadratic by
the theorem of cube and \lq{pointe de degree 2\rq \ in the sense
of [MB], p.12. Indeed $\delta(rs)=r^2\delta(s)$ for all $r\in \N$,
since ${\cal L}$ is symmetric. If $m$ annihilates ${\cal H}({\cal
L}^2)$ the values of $\delta$ are $m$-torsion (resp $2m$-torsion
if $2\vert m$) line bundles on $S$ by [MB], lemma 5.6. In our case
this applies for $m=4\nu$. If we extend the base to level $2m^2$,
then every $s\in {\cal H}({\cal L}^2)$ can be written in the form
$s=2my$. But $\delta(s)=\delta(2my)=4m^2\delta(y)$. Since
$\delta(y)$ is $4m^2$-torsion, we get $\delta(s)=0$. Hence a
$\delta$-marking exists after a suitable choice of level, i.e.
after an etale extension of our base scheme $S=A_{g-1,\kappa,N}$.

\bigskip\noindent
h) {\it The induction step}. Suppose $g\geq 3$. To show
$\Omega^1(X)=0$ it suffices to prove $u_*(\eta)^\chi=0$ for all
$\chi(t)=t^\nu$ with $\nu\geq 0$ and $\nu\not\equiv 0$ mod $p$ by
step e) and f). Indeed the Fourier-Jacobi expansion corresponds to
the expansion with respect to the parameter $q=q_{gg}$. So it is
enough to show $\Gamma({\cal A}_{g-1},{\cal E}\otimes {\cal L})$
for the separable invertible sheaves ${\cal L}={\cal L}_0^\nu$.
Using the exact sequence for ${\cal E}$ and the sequence $0\to
\Omega^1_{X} \to \Omega^1_{{\cal A}_{g-1}} \to \Omega_{A}\to 0 $
for $A={\cal A}_{g-1}$ over $S$
 it suffices to show for $X= A_{g-1,\kappa,N}$ the vanishing
$\Gamma(X,K)=0$, $\Gamma(X,\Omega_A\otimes K)=0$ and
$\Gamma(X,\Omega^1_X\otimes K)=0$. By $K\hookrightarrow {\cal
O}_X$ and by the Kodaira-Spencer isomorphism the last two
statements are reduced to the induction assumption
$\Gamma(X,\Omega_X^1)=0$ in genus $g-1$. This proves the induction
step, since obviously $\Gamma(X,K)=0$.

\bigskip\noindent
i) {\it Start of induction}. We now consider the case $g=2$. Here
the argument of step h) fails, since the space of 1-forms in genus
one usually is nontrivial. However in the situation as in h) the
argument implies that under the projection $r:{\cal E}\otimes{\cal
L}\to {\cal L}$
$$ r(u_*(\eta)^\chi) \ \in \ \Gamma({\cal A}_{g-1},{\cal L}) $$
vanishes. Indeed, this is a section in $\Gamma({\cal
A}_{g-1},{\cal L}) \cong \Gamma(X,\pi_*({\cal L}))$. Let us
consider its Fourier expansion. On can easily reduce to consider
the standard cusp at infinity by replacing $\eta$ without changing
the condition obtained by c). By the Koecher effect applied for
$\eta$ this expansion is holomorphic at the cusp. Hence
$r(u_*(\eta)^\chi)$ vanishes, since a power of it is a modular
form of negative weight. Now this implies that all $
f(q_{11},q_{12}) q_{22}^\nu \cdot \frac{dq_{22}}{q_{22}} $ terms
in the Fourier expansion of $\eta$ vanish for $\nu \not\equiv 0 $
mod $p$. In other words
$$ (**) \quad \quad
\fbox{$ {}^t \overline v \overline T \overline v \ \neq\ 0 \
\mbox{ for }\ \overline v\in \F_p^g \quad \Longrightarrow \quad
{}^t \overline v \overline a(T) \overline v\ =\ 0 $} \ .$$ In fact
by step b) it is enough to have this for ${}^t\overline v =
(0,1)$, since in general ${}^t v= (0,1) {}^tU$ for some unimodular
matrix $U$ congruent to $id$ mod $N$.

\bigskip\noindent
j) Suppose $det(T)\neq 0$ mod $p$. Then the reduction $\overline
T$ of $T$ mod $p$ is a nondegenerate binary quadratic form over
$\F_p$. Also $S=a(T)$ is a binary quadratic form, and it is easy
to see that (**) then implies $S=0$. Also $a(0)=0$ by step b),
hence $a(T)$ vanishes unless the reduction $\overline T$ of $T$
mod $p$ has rank 1 with coefficient $\nu = \overline T_{22}\neq 0$
in $\F_p\subseteq \kappa$. Then necessarily
$$ \overline T = \begin{pmatrix} x^2\nu & x\nu \cr
x\nu & \nu \cr
\end{pmatrix} = \begin{pmatrix} 1 & x \cr
0 & 1 \cr
\end{pmatrix}
\begin{pmatrix} 0 & 0 \cr 0 & \nu \cr
\end{pmatrix}
\begin{pmatrix} 1 & 0 \cr x & 1 \cr
\end{pmatrix}  $$
for some $x\in \F_p\subseteq \kappa$.

\bigskip\noindent
k) {\it Fourier-Jacobi expansion revisited}. We now compare the
Fourier expansion of $\eta\in \Omega^1(X)$ with its Fourier-Jacobi
expansion for $q=q_{22}$.  To unburden notation let us assume
level $N=1$ for simplicity. Then
$$ \eta = \sum_{\nu=0}^\infty \eta_\nu(q_{11},q_{12}) q_{22}^\nu \ ,$$
where $\eta_\nu = \sum A(\begin{smallmatrix} t_0 & t_1 \cr t_1 &
\nu \cr
\end{smallmatrix}) q_{11}^{t_0}q_{12}^{2t_1}$ with summation over
semidefinite matrices $(\begin{smallmatrix} t_0 & t_1 \cr t_1 &
\nu \cr
\end{smallmatrix})$ with $t_0\in \Z$ and $2t_1\in \Z$.
Using the transformation property b) for $U= (\begin{smallmatrix}
1 & 0 \cr g & 1 \cr
\end{smallmatrix})$ with $g\in \Z$ we can also write
$$ \eta_\nu(q_{11},q_{12}) = \sum_{a}\bigl(\sum_{t_0} \vartheta_a^{A(t_0)} \cdot q_{11}^{t_0 -
\frac{a^2}{4\nu}}\bigr) $$ where now $a$ runs over the finitely
many cosets $a\in \frac{1}{2}\Z^{g-1}/\nu \Z^{g-1}$, and where
$t_0$ runs over all integers $t_0\geq \frac{a^2}{4\nu}$. Notice
that $A(t_0) = (\begin{smallmatrix} t_0 & a \cr a & \nu \cr
\end{smallmatrix})$ depends on $t_0$. Here we used
for an integral matrix $2\times 2$-matrix $A=(\begin{smallmatrix}
a_0 & a_1 \cr a_1 & a_2 \cr
\end{smallmatrix})$ (or its reduction
mod $p$) the notation
$$ \vartheta_a^A = \sum_{g\in \Z} \begin{pmatrix} 1 & g \cr 0 & 1 \cr
\end{pmatrix} A\begin{pmatrix} 1 & 0 \cr g & 1 \cr
\end{pmatrix} \cdot q_{11}^{\nu(g+\frac{a}{\nu})^2} q_{12}^{2\nu(g+\frac{a}{\nu})}
$$ $$ = \sum_{g\in \Z} \begin{pmatrix} a_0+2ga_1 + g^2a_2 & a_1+ga_2
\cr
 a_1+ga_2 & a_2 \cr
\end{pmatrix} \cdot q_{11}^{\nu(g+\frac{a}{\nu})^2} q_{12}^{2\nu(g+\frac{a}{\nu})}\ .$$
$$ = \begin{pmatrix} a_0 & a_1
\cr
 a_1 & a_2 \cr
\end{pmatrix} f_a \ +\
\begin{pmatrix} 2a_1 & a_2
\cr
 a_2 & 0 \cr
\end{pmatrix} f_a^{(1)}
\ +\
\begin{pmatrix} a_2 & 0
\cr
 0 & 0 \cr
\end{pmatrix} f_a^{(2)} $$
for the derivatives $f_a^{(i)} = \sum_{g\in \Z} g^i
q_{11}^{\nu(g+\frac{a}{\nu})^2}q_{12}^{2\nu(g+\frac{a}{\nu})}$ of
the theta functions $f_a(q_{11},q_{12})$. Summation over $t_0$
therefore expresses $\eta_\nu(q_{11},q_{12})$ as a sum over $a$ of
terms of the form
$$ f_a \cdot \bigl(
g_{0,a} \frac{dq_{11}}{q_{11}} + 2 g_{1,a}\frac{dq_{12}}{q_{12}} +
g_{2,a}\frac{dq_{22}}{q_{22}} \bigr)\quad + $$
$$
 f_a^{(1)} \cdot \bigl( 2 g_{1,a}\frac{dq_{11}}{q_{11}} +
2g_{2,a}\frac{dq_{12}}{q_{12}} \bigr) \ \ +\ \ f_a^{(2)} \cdot
g_{2,a}\frac{dq_{11}}{q_{11}}
 $$
for power series $g_{i,a}=g_{i,a}(q_{11})$ in the variable
$q_{11}$. Linear independency of the theta functions $f_a$
(Mumford's theory) and step i) imply
$$ g_{2,a}(q_{11}) = 0 $$
by looking at the coefficient at
$\frac{dq}{q}=\frac{dq_{22}}{q_{22}}$. Hence the terms simplify to
$$ f_a \cdot \bigl(
g_{0,a} \frac{dq_{11}}{q_{11}} + 2 g_{1,a}\frac{dq_{12}}{q_{12}}
\bigr)\quad + \quad
 f_a^{(1)} \cdot  2 g_{1,a}\frac{dq_{11}}{q_{11}} \ .
 $$
The coefficients $g_{1,a}$ are sections in
$\Gamma(A_{1,\kappa,n},\Omega_{A_{1,\kappa,n}} \otimes K^{-1})$
and the coefficients $g_{0,a}$ are sections in
$\Gamma(A_{1,\kappa,n},\Omega^1_{A_{1,\kappa,n}} \otimes K^{-1})$
for suitably levels $n$ as in step h). For the choice of $n$ see
step f) and g). Morally speaking $g_{1,a},g_{0,a}$ are modular
forms of genus one of level $n$ and weight $\frac{1}{2}$
respectively $\frac{3}{2}$, if this is suitably defined via theta
line bundles. By the way they arise it is easy to see, that they
are holomorphic at the cusps.

\bigskip\noindent
l) {\it $p$-singularity}. Notice the Fourier coefficients of
$\vartheta^{A(t_0)}q_{22}^\nu$ are all of the form
$$ T'= \begin{pmatrix} \nu y^2 & \nu y \cr
\nu y & \nu \cr \end{pmatrix} $$ where $y \in \frac{1}{M}\Z$ for
some integer $M$ with $(M,p)=1$. Notice $M$ also takes care of the
level $N$, which was suppressed so far. Looking at the Fourier
coefficients of $\frac{dq_{12}}{q_{12}}$ in the $q$-expansion of
$\eta$ only those $T$ contribute with $T = T' +
(\begin{smallmatrix} l & 0 \cr 0 & 0 \cr \end{smallmatrix})$,
where $l$ comes from the Fourier expansion of the section
$g_{1,a}(q_{11})$. Reading this modulo $p$ by step j) gives
$$ \begin{pmatrix} \nu x^2 & \nu x \cr
\nu x & \nu \cr \end{pmatrix} \ \equiv \ \begin{pmatrix} \nu y^2 &
\nu y \cr \nu y & \nu \cr \end{pmatrix} + \begin{pmatrix} l & 0
\cr  0 & 0 \cr \end{pmatrix} $$ modulo $p$, since $\nu\not\equiv
0$ mod $p$. Hence $x$ and $y$ are congruent modulo $p$, and
therefore
$$ l \equiv 0 \mod p \ .$$
This implies that the $q_{11}$-expansions of the sections
$g_{1,a}$ are $p$-singular for all $a$. Using this consider the
$\frac{dq_{11}}{g_{11}}$ term. From the linear independency of the
$f_a$ then the same also follows  for the $g_{0,a}$. In other
words
$$ g_{i,a}(q_{11}) = \sum_{n\geq 0} c_{i,a,n}\cdot q_{11}^{n/M} \quad , \quad
c_{i,a,n}=0 \mbox{ if } \ p\nmid n \ .$$

\bigskip\noindent
m) {\it Vanishing of the $g_{i,a}$}.
 Since $(g_{i,a})^2$ are modular
forms of weight 1 resp. 3 of genus one with suitable level $n$ in the sense of Katz [K2] by
step k), the $p$-singularity shown in step l) contradicts corollary (2) of the theorem of loc.
cit, p.55 unless $(g_{i,a})^2=0$ since 1 resp. 3 is $<p-1$ by our assumption $p\neq 2,3$. This
theorem of Katz is a version in the elliptic case with level structure of part 2) of our second
corollary above. We conclude $g_{i,a}=0$ for all $i$, hence $\eta=0$. This proves the vanishing
theorem for $p\neq 2,3$. For $p=3$ we only have $g_{1,a}=0$. But this implies, that all Fourier
coefficients $A(T)$ for matrices $T$ with $T_{22} \not\equiv 0 $ mod $3$ have coefficients
$A(T)_{ij}=0$ unless $i=1,j=1$. But this implies $A(T)=0$ for all such Fourier coefficients by
the transformation property (*). By j) it is enough to suppose $\overline T =
(\begin{smallmatrix} 0 & 0 \cr 0 & 1\cr\end{smallmatrix})$. Then otherwise up to a constant
$A(T)= (\begin{smallmatrix} 1 & 0 \cr 0 & 0\cr\end{smallmatrix})$. Now apply $U$ with
$\overline U = (\begin{smallmatrix} 0 & -1 \cr 1 & 1\cr\end{smallmatrix})$. Then we get
$\overline T_1 = (\begin{smallmatrix} 1 & 1 \cr 1 & 1\cr\end{smallmatrix})$ and $A(T_1) =
(\begin{smallmatrix} 0 & 0 \cr 0 & 1\cr\end{smallmatrix})$ from (*). A contradiction to (**).
This proves the theorem for $p=3$. \hfill $\square$

\bigskip\noindent
{\it Proof of the unique factorization theorem}. We may suppose
$g\geq 3$, since for $g=2$ the assertion follows from [I]. Let
$R_N=M_{g,\kappa,N}$ denote the ungraded ring of modular form of
level $g$ with level-$N$-structure. We  have to show that the
class group $Cl(R_N)$ of the Krull domain $R_N=M_{g,\kappa,N}$
vanishes for level $N=1$. For suitably high level $N$ (divisible
by 8) there exists a projective embedding $\phi: A_{g,\kappa,N}^*
\hookrightarrow {\mathbb P}$ of the minimal compactification
$X^*=A_{g,\kappa,N}^*$ of $X=A_{g,\kappa,N}$ into projective space
by the theory of theta functions, such that $j^*({\cal O}_{\mathbb
P}(1)) = \omega$. See [CF], p.157-159  and [CF], V.5. (It is
preferable to use theta-level-structures $(N,2N)$, but we skip a
discussion of this).

\bigskip\noindent By definition $Spec(R_N) =C(X^*)$ is the affine cone over
the projective subvariety $j:X^*=A_{g,\kappa,N}^*\hookrightarrow
\mathbb P$, since it is the ungraded ring of the corresponding
graded coordinate ring of $X^*$. Therefore $Cl(R_N) = Cl(C(X^*))$
is obtained from $Cl(X^*)$ by
$$   Cl(C(X^*)) = Cl(X^*)\big/\Z\cdot H $$
where $H$ is the class corresponding to a hyperplane section $ H$,
which  is a section of $j^*({\cal O}_{\mathbb P}(1)) = \omega$.
The cone $C(X^*)$ is normal and regular outside the isolated cusp
and the singular lines over $X^*_s=X^*\setminus X^*_{reg}$. Since
$g\geq 3$, the locus $X_s^*$ has codimension $\geq 2$. Hence there
exist canonical $G_N$-equivariant isomorphisms $Cl(C(X^*)) \cong
Pic(C(X^*))$ and $Cl(X^*) = Pic(X^*)$ and $Cl(X)\cong Pic(X)$. The
complement of $X_s = X\setminus X_{reg}$ in $X_s^*$ also has
codimension 2. Hence also $Cl(X) \cong Cl(X^*)$, and therefore
$Pic(X^*)=Pic(X)$. The class of $H$ therefore corresponds to the
class of $\omega$ in $Pic(X)$.
 Hilbert theorem 90
implies
$$  Cl(R_1) \hookrightarrow Cl(R_N)^{G_N}  \ .$$
The composed map $\delta: Cl(R_N)^{G_N} \to H^1(G_N,Quot(R_N)^*/\kappa^*) \to
H^2(G_N,\kappa^*)$ is trivial on the image of $Cl(R_1)$. A class $x$ in $Cl(R_1)$ defines a
class $y$ in $Cl(R_N)$ with vanishing obstruction in $H^2(G_N,\kappa^*)$ (notice for $2\vert N$
this group may be nontrivial). Via the isomorphisms above there exists a corresponding class
$y^*$ in $Pic(X)/\Z\cdot c_1(\omega)\cong Pic(X^*)/\Z\cdot c_1(\omega) \cong Cl(X^*)/\Z\cdot
j^*({\cal O}(1)) \cong Cl(C(X^*))$, whose obstruction in $H^2(G_N,\kappa^*)$ vanishes in the
sense of step 1 of the proof of the theorem. Hence by descend any representative $y\in
Pic(A_{g,\kappa,N})$ of $y^*$ is the pullback of a class from $Pic(A_{g,\kappa})$. However
$Pic(A_{g,\kappa}) = \Z \cdot c_1(\omega)$ as shown in the theorem for $p>p_0$. Therefore
$y^*=0$, hence $y=0$ and $x=0$. This proves $Cl(R_1)=0$ for $p>p_0$, hence the unique
factorization theorem. \hfill$\square$

\newpage

\bigskip\noindent
\centerline{\bf References}


\bigskip\noindent
[C] Chai C.L., Compactification of Siegel moduli schemes, London
Math. Soc. Lecture Notes Series 107, Cambridge university press,
(1985)

\bigskip\noindent
[CF] Chai C.L.-Faltings G., Degeneration of Abelian Varieties,
Ergebnisse der Mathematik und ihrer Grenzgebiete, vol. 22,
Springer Verlag (1990)

\bigskip\noindent
[D] Demazure M., Lectures on $p$-Divisible Groups, SLN 302,
Springer Verlag (1972)

\bigskip\noindent
[I] Igusa J.I., On the ring of modular forms of degree two over
$\Z$, Jorn. of Math. 101, 1-3 (1979), 149 - 183

\bigskip\noindent
[FD] Fiederowisz Z.-Priddy S., Homology of Classical Groups over a
Finite Field, in Algebraic K-theory (Evanston 1976), edited by
M.R.Stein, SLN 551, Springer Verlag, p. 269 - 282

\bigskip\noindent
[F] Freitag E., Die Irreduzibilit\"{a}t der Schottkyrelation
(Bemerkung zu einem Satz von J.Igusa), Archiv Math. vol. 40
(1983), 255 - 259

\bigskip\noindent
[K] Katz N.M. $p$-adic properties of modular schemes and modular
forms, Modular Functions of One Variable III, SLN vol. 350,
Springer Verlag

\bigskip\noindent
[K2] Katz N.M. A result of modular forms in characteristic $p$,
Modular Functions of One Variable V, SLN vol. 601, Springer Verlag

\bigskip\noindent
[MB] Moret-Bailly, L., Pinceaux de varietes abeliennes, Asterisque
129, (1985)

\bigskip\noindent
[M] Mumford D., On the equations defining abelian varieties I,
Invent. Math. 1 (1966), 287 - 354129 (1985)

\bigskip\noindent
[M2] Mumford D., On the equations defining abelian varieties II,
Invent. Math. 3 (1967), 75 - 135

\bigskip\noindent
[S] J.P.Serre, Congruences et formes modulaires, Seminaire
Bourbaki 1971/72, exp. 416, SLN 317


\bigskip\noindent
[Sw] Swinnerton-Dyer H.P.F., On $l$-adic representations and
congruences for coefficients of modular forms, Modular functions
in one Variable III, SLN vol. 350, Springer Verlag

\end{document}